\documentclass{article} 
\usepackage{mathrsfs}
\usepackage{amsmath}
\usepackage{amssymb}
\usepackage{amsfonts}
\usepackage{amsmath}
\usepackage[numbers]{natbib}
\usepackage[left=1in, right=1in, top=1in, bottom=1in]{geometry}
\setcitestyle{open={},close={}}

\allowdisplaybreaks
\begin{document}
\author{{\bf X. Shi$^{1}$, V.D. R\u{a}dulescu$^{2}$, D.D. Repov\v{s}$^{3}$, and Q. Zhang$^{4}$} \bigskip\\ 
$^{1}${\small {Zhengzhou University of Light Industry, Zhengzhou, Henan, China. \texttt{aryang123@163.com}}}\\
$^{2}${\small {AGH University of Science and Technology, Krak\'ow, Poland.   \texttt{vicentiu.radulescu@math.cnrs.fr}}}\\
$^{3}${\small {University of Ljubljana, Ljubljana, Slovenia.  \texttt{dusan.repovs@guest.arnes.si}}}\\
$^{4}${\small {University of Georgia, Athens, GA, 30602, USA.  \texttt{zhangqihu@yahoo.com}}}}
\title{\bf Multiple solutions of double phase variational problems
with variable exponent}
\date{}
\maketitle
\begin{abstract}
This paper deals with the existence of multiple solutions for the
quasilinear equation $-\mathrm{div}\,\mathbf{A}(x,\nabla u)+\left\vert
u\right\vert ^{\alpha (x)-2}u=f(x,u)$ in $
\mathbb{R}
^{N}$, which involves a general variable exponent elliptic operator $\mathbf{
A}$ in divergence form. The problem corresponds to double phase anisotropic
phenomena, in the sense that the differential operator has behaviors like $
\left\vert \xi \right\vert ^{q(x)-2}\xi $ for small $\left\vert \xi
\right\vert $ and like $\left\vert \xi \right\vert ^{p(x)-2}\xi $ for large $
\left\vert \xi \right\vert $, where $1<\alpha (\cdot )\leq p(\cdot )<q(\cdot
)<N$. Our aim is to approach variationally the problem by using the tools of
critical points theory in generalized Orlicz-Sobolev spaces with variable
exponent. Our results extend the previous works Azzollini, d'Avenia, and  Pomponio (2014) and  Chorfi and  R\u{a}dulescu (2016),
  from the case when
exponents $p$ and $q$ are constant, to the case when $p(\cdot )$ and $%
q(\cdot )$ are functions. We also substantially weaken some of their
hypotheses overcome the lack of compactness by using
the weighting method.\\


\end{abstract}
\section{Introduction}

In this paper, we deal with the following variable exponent elliptic
equation
\begin{equation*}
-\mathrm{div}\,\mathbf{A}(x,\nabla u)+\left\vert u\right\vert ^{\alpha
(x)-2}u=f(x,u):=\lambda a(x)\left\vert u\right\vert ^{\delta (x)-2}u+\mu
w(x)g(x,u),\text{ \ }\left. {}\right. \eqno(\mathcal{E})
\end{equation*}%
where $\lambda >0$ and $\mu \geq 0$ are parameters, $\mathbf{A:}$ $%
\mathbb{R}
^{N}\times
\mathbb{R}
^{N}\rightarrow
\mathbb{R}
^{N}$ admits a potential $\mathscr{A}$, with respect to its second variable $%
\xi $, satisfying the following assumption:

\smallskip ($\mathcal{A}_{1}$) the potential $\mathscr{A}=\mathscr{A}(x,\xi
) $ is a continuous function in $%
\mathbb{R}
^{N}\times
\mathbb{R}
^{N}$, with continuous derivative with respect to $\xi $, $\mathbf{A}%
=\partial _{\xi }\mathscr{A}(x,\xi )$, and verifies:

(i) $\mathscr{A}(x,0)=0$ and $\mathscr{A}(x,\xi )=\mathscr{A}(x,-\xi )$, for
all $(x,\xi )\in
\mathbb{R}
^{N}\times
\mathbb{R}
^{N}$;

(ii) $\mathscr{A}(x,\cdot )$ is strictly convex in $%
\mathbb{R}
^{N}$ for all $x\in $ $%
\mathbb{R}
^{N}$;

(iii) there exist positive constants $C_{1},\,C_{2}$ and variable exponents $%
p$ and $q$ such that for all $(x,\xi )\in
\mathbb{R}
^{N}\times
\mathbb{R}
^{N}$%
\begin{equation}
\left.
\begin{array}{c}
C_{1}\left\vert \xi \right\vert ^{p(x)},\text{ if }\left\vert \xi
\right\vert >>1 \\
C_{1}\left\vert \xi \right\vert ^{q(x)},\text{ if }\left\vert \xi
\right\vert <<1%
\end{array}%
\right\} \leq \mathbf{A}(x,\xi )\cdot \xi \quad\text{and}\quad\left\vert
\mathbf{A}(x,\xi )\right\vert \leq \left\{
\begin{array}{c}
C_{2}\left\vert \xi \right\vert ^{p(x)-1},\text{ if }\left\vert \xi
\right\vert >>1 \\
C_{2}\left\vert \xi \right\vert ^{q(x)-1},\text{ if }\left\vert \xi
\right\vert <<1%
\end{array}%
;\right.  \label{5.1}
\end{equation}

(iv) $1<<p(\cdot )<<q(\cdot )<<\min \{N,p^{\ast }(\cdot )\}$, and $p(\cdot
), $ $q(\cdot )$ are Lipschitz continuous in $%
\mathbb{R}
^{N}$;

(v) $\mathbf{A(}x,\xi )\cdot \xi \leq s(x)\mathscr{A}(x,\xi )$ for any $%
(x,\xi )\in
\mathbb{R}
^{2N}$, where $s$ is Lipschitz continuous and satisfies $q(\cdot )\leq
s(\cdot )<<p^{\ast }(\cdot )$;

\smallskip ($\mathcal{A}_{2}$) $\mathscr{A}$ is uniformly convex, that is,
for any $\varepsilon \in (0,1)$, there exists $\delta (\varepsilon )\in
(0,1) $ such that $\left\vert u-v\right\vert \leq \varepsilon \max
\{\left\vert u\right\vert ,\left\vert v\right\vert \}$ or $\mathscr{A}(x,%
\frac{u+v}{2})\leq \frac{1}{2}(1-\delta (\varepsilon ))(\mathscr{A}(x,u)+%
\mathscr{A}(x,v)) $ for any $x,u,v\in
\mathbb{R}
^{N}$. \newline

In this paper, for any $v:%
\mathbb{R}
^{N}\rightarrow
\mathbb{R}
$, we denote
\begin{equation*}
v^{+}=\underset{x\in
\mathbb{R}
^{N}}{\text{ess}\sup }\,v(x)\text{, }v^{-}=\underset{x\in
\mathbb{R}
^{N}}{\text{ess}\inf }\,v(x)\text{,}
\end{equation*}%
and we denote by $v_{1}<<v_{2}$ the fact that
\begin{equation*}
\underset{x\in
\mathbb{R}
^{N}}{\text{ess}}\inf (v_{2}(x)-v_{1}(x))>0\text{.}
\end{equation*}

\medskip
\textbf{Remark 1.1}
 A typical case of $\mathbf{A}$ is
\begin{equation*}
 \mathbf{A}(x,\nabla
u)\mathbf{=}\left\{
\begin{array}{c}
\left\vert \nabla u\right\vert ^{p(x)-2}\nabla u\text{, if }\left\vert
\nabla u\right\vert >1 \\
\left\vert \nabla u\right\vert ^{q(x)-2}\nabla u\text{, if }\left\vert
\nabla u\right\vert \leq 1
\end{array}
\right.   
\end{equation*}

Then
\begin{equation*}
-\mathrm{div}\, \mathbf{A}(x,\nabla u)=\left\{
\begin{array}{c}
-\mathrm{div}\, (\left\vert \nabla u\right\vert ^{p(x)-2}\nabla u)\text{, if
}\left\vert \nabla u\right\vert >1 \\
-\mathrm{div}\, (\left\vert \nabla u\right\vert ^{q(x)-2}\nabla u)\text{, if
}\left\vert \nabla u\right\vert \leq 1
\end{array}
\right.
\end{equation*}

and
\begin{equation*}
\mathscr{A}(x,\xi )=\left\{
\begin{array}{c}
\frac{1}{p(x)}\left\vert \xi \right\vert ^{p(x)}+\frac{1}{q(x)}-\frac{1}{p(x)
}\text{, if }\left\vert \xi \right\vert >1 \\
\frac{1}{q(x)}\left\vert \xi \right\vert ^{q(x)}\text{, if }\left\vert \xi
\right\vert \leq 1.
\end{array}
\right.
\end{equation*}
From Lemma A.2, given in Appendix A of [\cite{w12}], it is clear that this
typical potential $\mathscr{A}$ satisfies ($\mathcal{A}_{1}$) and ($\mathcal{
A}_{2}$), $1<p^{-}\leq p^{+}<N$ and $1<q^{-}\leq q^{+}<N$.\\

It is well known that a main difficulty in studying the elliptic equations
in $\mathbb{R}^{N}$ is the lack of compactness. To overcome this difficulty, many methods
can be used. One type of methods is that under some additional conditions
there holds the required compact imbedding theorem, for example, the
symmetry method, the coercive coefficient method and the weighting method.
In [\cite{w2}], the authors consider the equation $(\mathcal{E})$ with
constant exponent by the symmetry method to rebuild the required compact
imbedding theorem. In [\cite{w12}], in order to rebuild the required compact
imbedding theorem, the authors consider the equation $(\mathcal{E})$ with
coercive coefficient $V(x)$ of $\left\vert u\right\vert ^{\alpha (x)-2}u$,
namely, $V(x)\rightarrow +\infty $ as $|x|\rightarrow \infty $. In this
paper, we will apply the weighting method, namely if the coefficients $w$ and $%
a $ satisfy some integrable conditions, then we can rebuild the required
compact embedding theorem.

We also make the following assumptions:

\smallskip ($\mathcal{H}_{f}^{1}$) $g:
\mathbb{R}
^{N}\times
\mathbb{R}
\rightarrow
\mathbb{R}
$ satisfies the Carath\'eodory condition, $0\leq g(x,u)u=o(\left\vert
u\right\vert ^{\alpha (x)})$ as $u\rightarrow 0$, and $\left\vert
g(x,u)\right\vert \leq C(1+\left\vert u\right\vert ^{\gamma (x)-1})$, where $
\gamma (\cdot )$ is Lipschitz continuous and $\alpha \leq \gamma (\cdot
)<<p^{\ast }(\cdot )$.

\smallskip ($\mathcal{H}_{f}^{2}$) There exists a constant $\theta >s^{+}$
such that
\begin{equation*}
0<G(x,t)\leq \frac{1}{\theta }tg(x,t),\quad \forall t\in
\mathbb{R}
\backslash \{0\},\ \forall x\in \mathbb{R}^{N},
\end{equation*}%
where $G(x,t)=\int_{0}^{t}g(x,s)ds$, and $s(\cdot )$ is defined in ($
\mathcal{A}_{1}$)-(v).

\smallskip ($\mathcal{H}_{f}^{3}$) $g(x,-u)=-g(x,u)$.

\smallskip ($\mathcal{H}_{w}$) $w\in L^{r(\cdot )}(
\mathbb{R}
^{N})$, $w>0$ a.e. in $\mathbb{R}
^{N}$, $1<<r(x)<<\infty $, and
\begin{equation*}
r^{\prime }(x)\leq \frac{p^{\ast }(x)}{\gamma (x)},\quad \forall x\in
\mathbb{R}
^{N},
\end{equation*}
where $r^{\prime }(x)$ is the conjugate function of $r(x)$, namely $\frac{1}{
r(x)}+\frac{1}{r^{\prime }(x)}=1$, and
\begin{equation*}
p^{\ast }(x)=\left\{
\begin{array}{l}
Np(x)/(N-p(x))\text{ , }\quad\mbox{if}\ p(x)<N, \\
\infty \text{ , }\quad\mbox{if}\ p(x)\geq N.
\end{array}
\right.
\end{equation*}

($\mathcal{H}_{a}$) $1<<\delta (x)\leq \delta ^{+}<\alpha ^{-}$, $a\in
L^{r_{\ast }(\cdot )}(
\mathbb{R}
^{N})$, $a>0$ a.e. in $
\mathbb{R}
^{N}$, $1<<r_{\ast }(x)<<\infty $, and
\begin{equation*}
\alpha (x)\leq \frac{r_{\ast }(x)}{r_{\ast }(x)-1}\delta (x)\leq p^{\ast
}(x),\quad\forall x\in
\mathbb{R}^{N}.
\end{equation*}

\smallskip This paper generalizes some results contained in [\cite{w2,w7}] to the case of partial differential equations with variable
exponent. If $p(\cdot )\equiv p$, $q(\cdot )\equiv q$ and $\alpha (\cdot)\equiv \alpha $ are constants, then $(\mathcal{E})$ becomes the usual
constant exponent differential equation in divergence form discussed in [\cite{w7}]. But if either $p(\cdot )$ or $q(\cdot )$ are nonconstant
functions, then $(\mathcal{E})$ has a more complicated structure, due to its
non-homogeneities and to the presence of several nonlinear terms.

This paper is motivated by double phase nonlinear problems with variational
structure, which have been introduced by Marcellini [\cite{marce1}] and
developed by Mingione \textit{et al.} [\cite{mingi1}, \cite{mingi2}] in the
framework of nonhomogeneous problems driven by a differential operator with
variable growths described by nonconstant functions $p(x)$ and $q(x)$. In
the case of two different materials that involve power hardening exponents $p(\cdot )$ and $q(\cdot )$, the differential operator $\mathrm{div}\,\mathbf{A}(x,\nabla u)$ describes the geometry of a composite of these two
materials. Cf. hypothesis \eqref{5.1}, the $p(\cdot )$-material is present
if $|\xi |>>1$. In the contrary case, the $q(\cdot )$-material is the only
one describing the composite.

In recent years, the study of differential equations and variational
problems with variable exponent growth conditions have been an interesting
topic, which have backgrounds in image processing, nonlinear
electrorheological fluids and elastic mechanics etc. We refer the readers to
[\cite{e1, e2, e18, j6, radrep, e3, e4}] and the references therein for more background of
applications. There are many reference papers related to the the study of
variational problems with variable exponent growth conditions, far from
being complete, we refer readers to [\cite{jj1,e7,j22,j8}],[\cite{e8}-\cite{14}], [\cite{34}-\cite{19a}], [\cite{radnla,repovs}], [\cite{e26}-\cite{j7}], [\cite{e29,w12}].

\smallskip
Our main results can be stated as follows.

\medskip
\textbf{Theorem 1.1}. Assume that $1<<\alpha \leq p<<q<<\min \{N,p^{\ast
}\}$, $1<<\alpha (\cdot )<<p^{\ast }(\cdot )\frac{q^{\prime }(\cdot )}{p^{\prime }(\cdot )}$, $\mu >0$, $\lambda $ is small enough, and hypotheses ($\mathcal{A}%
_{1}$)-($\mathcal{A}_{2}$), ($\mathcal{H}_{f}^{1}$)-($\mathcal{H}_{f}^{2}$),
($\mathcal{H}_{w}$) and ($\mathcal{H}_{a}$) hold. Then problem $(\mathcal{E})$
has two pairs of nontrivial non-negative and nonpositive solution.\newline

\textbf{Theorem 1.2}. Assume that $1<<\alpha \leq p<<q<<\min \{N,p^{\ast
}\}$, $1<<\alpha (\cdot )<<p^{\ast }(\cdot )\frac{q^{\prime }(\cdot )}{p^{\prime }(\cdot )}$, $\mu >0$, and hypotheses ($\mathcal{A}_{1}$)-($\mathcal{A}_{2}$), ($%
\mathcal{H}_{f}^{1}$)-($\mathcal{H}_{f}^{3}$), ($\mathcal{H}_{w}$) and ($\mathcal{H}_{a}$) hold. Then problem $(\mathcal{E})$ has infinitely many
nontrivial solutions with energy tending to $+\infty $.\newline

\textbf{Theorem 1.3}. Assume that $1<<\alpha \leq p<<q<<\min \{N,p^{\ast
}\}$, $1<<\alpha (\cdot )<<p^{\ast }(\cdot )\frac{q^{\prime }(\cdot )}{p^{\prime }(\cdot )}$, $\mu =0$, and hypotheses ($\mathcal{A}_{1}$)-($\mathcal{A}_{2}$), ($%
\mathcal{H}_{w}$) and ($\mathcal{H}_{a}$) hold. Then problem $(\mathcal{E})$ has
infinitely many nontrivial solutions with negative energy tending to $0$.\newline

This paper is divided into five sections; Section 2 contains some properties
of function spaces with variable exponent; Section 3 includes several basic
properties of Orlicz-Sobolev spaces; in Section 4 we establish some
qualitative properties of the operators involved in our analysis; Section 5
give the proofs of Theorems 1.1-1.3, respectively.
We refer to [\cite{ciarlet}] for the basic analytic tools used in this paper.

\section{Variable Exponent Spaces Theory}

Nonlinear problems with non-homogeneous structure are motivated by numerous
models in the applied sciences that are driven by partial differential
equations with one or more variable exponents. In some circumstances, the
standard analysis based on the theory of usual Lebesgue and Sobolev function
spaces, $L^{p}$ and $W^{1,p}$, is not appropriate in the framework of
materials that involve non-homogeneities. For instance, both
electrorheological \textquotedblleft smart fluids" fluids and phenomena
arising in image processing are described in a correct way by nonlinear
models in which the exponent $p$ is not necessarily constant. The variable
exponent describe the geometry of a material which is allowed to change its
hardening exponent according to the point. This leads to the analysis of
variable exponents Lebesgue and Sobolev function spaces (denoted by $%
L^{p(\cdot )}$ and $W^{1,p(\cdot )}$), where $p$ is a real-valued
(non-constant) function.

Throughout this paper, the letters $c,c_{i},C,C_{i}$, $i=1,2,\ldots $ denote
positive constants which may vary from line to line but are independent of
the terms which will take part in any limit process.

In order to discuss the problem $(\mathcal{E})$, we need some theories on
variable exponent Lebesgue spaces and Sobolev spaces. In the following, we
will give some properties of these variable exponent spaces. Let $\Omega
\subset
\mathbb{R}^{N}$ be an open domain. Let $S(\Omega )$ be the set of all measurable real
valued functions defined on $\Omega $. Let
\begin{equation*}
C_{+}(\overline{\Omega })=\left\{ v\left\vert v\in C(\overline{\Omega })%
\text{, }v(x)>1\text{ for }x\in \overline{\Omega }\right. \right\} ,
\end{equation*}%
\begin{equation*}
L^{p(\cdot )}(\Omega )=\left\{ u\in S(\Omega )\mid \int_{\Omega }\left\vert
u(x)\right\vert ^{p(x)}dx<\infty \right\} .
\end{equation*}

The function space $L^{p(\cdot )}(\Omega )$ is equipped with the Luxemburg
norm%
\begin{equation*}
\left\vert u\right\vert _{L^{p(\cdot )}(\Omega )}=\inf \left\{ \lambda
>0\left\vert \int_{\Omega }\left\vert \frac{u(x)}{\lambda }\right\vert
^{p(x)}dx\leq 1\right. \right\} .
\end{equation*}

Then ($L^{p(\cdot )}(\Omega )$, $\left\vert \cdot \right\vert _{L^{p(\cdot
)}(\Omega )}$) becomes a Banach space, we call it variable exponent Lebesgue
space.\newline

If $\Omega =%
\mathbb{R}
^{N}$, we simply denote ($L^{p(\cdot )}(%
\mathbb{R}
^{N}),\left\vert \cdot \right\vert _{L^{p(\cdot )}(%
\mathbb{R}
^{N})}$) as ($L^{p(\cdot )},\left\vert \cdot \right\vert _{L^{p(\cdot )}}$).%
\newline

\textbf{Proposition 2.1. }(see [\cite{7}, Theorem 1.15]) The space $%
(L^{p(\cdot )}(\Omega ),\left\vert u\right\vert _{L^{p(\cdot )}(\Omega )})$
is a separable, uniformly convex Banach space, and its conjugate space is $%
L^{p^{\prime }(\cdot )}(\Omega )$, where $\frac{1}{p(x)}+\frac{1}{p^{\prime
}(x)}=1$. For any $u\in L^{p(\cdot )}(\Omega )$ and $v\in L^{p^{\prime
}(\cdot )}(\Omega )$, we have the following H\"older inequality
\begin{equation*}
\left\vert \int_{\Omega }uvdx\right\vert \leq \left(\frac{1}{p^{-}}+\frac{1}{%
p^{\prime -}}\right)\left\vert u\right\vert _{L^{p(\cdot )}(\Omega )}\left\vert
v\right\vert _{L^{p^{\prime }(\cdot )}(\Omega )}.
\end{equation*}

\textbf{Proposition 2.2. }(see [\cite{7}, Theorem 1.16]) If $f:$ $\Omega
\times
\mathbb{R}
\rightarrow
\mathbb{R}
$ is a Carath\'{e}odory function and satisfies%
\begin{equation*}
\left\vert f(x,s)\right\vert \leq d(x)+b\left\vert s\right\vert
^{p_{1}(x)/p_{2}(x)}\text{ for any }x\in \Omega ,\ s\in
\mathbb{R}
,
\end{equation*}%
where $p_{1}$, $p_{2}\in C_{+}(\overline{\Omega })$, $d(x)\in L^{p_{2}(\cdot
)}(\Omega )$, $d(x)\geq 0$, $b\geq 0$, then the Nemytsky operator from $%
L^{p_{1}(\cdot )}(\Omega )$ to $L^{p_{2}(\cdot )}(\Omega )$ defined by $%
(N_{f}u)(x)=f(x,u(x))$ is a continuous and bounded operator. \newline

\textbf{Proposition 2.3. }(see [\cite{7}, Theorem 1.3]) If we denote
\begin{equation*}
\rho _{p(\cdot )}(u)=\int_{\Omega }\left\vert u\right\vert ^{p(x)}dx\text{, }%
\forall u\in L^{p(\cdot )}(\Omega ),
\end{equation*}%
then

i) $\left\vert u\right\vert _{L^{p(\cdot )}(\Omega
)}<1(=1;>1)\Longleftrightarrow \rho _{p(\cdot )}(u)<1(=1;>1);$

ii) $\left\vert u\right\vert _{L^{p(\cdot )}(\Omega )}>1\Longrightarrow
\left\vert u\right\vert _{L^{p(\cdot )}(\Omega )}^{p^{-}}\leq \rho _{p(\cdot
)}(u)\leq \left\vert u\right\vert _{L^{p(\cdot )}(\Omega )}^{p^{+}};$

$\left\vert u\right\vert _{L^{p(\cdot )}(\Omega )}<1\Longrightarrow
\left\vert u\right\vert _{L^{p(\cdot )}(\Omega )}^{p^{-}}\geq \rho _{p(\cdot
)}(u)\geq \left\vert u\right\vert _{L^{p(\cdot )}(\Omega )}^{p^{+}};$

iii) $\left\vert u\right\vert _{L^{p(\cdot )}(\Omega )}\rightarrow \infty
\Longleftrightarrow \rho _{p(\cdot )}(u)\rightarrow \infty .$\newline

\textbf{Proposition 2.4. }(see [\cite{7}, Theorem 1.4]) If $u$, $u_{n}\in
L^{p(\cdot )}(\Omega )$, $n=1,2,\cdots ,$ then the following statements are
equivalent:

1) $\underset{k\rightarrow \infty }{\lim }$ $\left\vert u_{k}-u\right\vert
_{L^{p(\cdot )}(\Omega )}=0;$

2) $\underset{k\rightarrow \infty }{\lim }$ $\rho _{p(\cdot )}\left(
u_{k}-u\right) =0;$

3) $u_{k}\rightarrow u$ in measure in $\Omega $ and $\underset{k\rightarrow
\infty }{\lim }$ $\rho _{p(\cdot )}\left( u_{k}\right) =\rho _{p(\cdot
)}\left( u\right) $.\newline

The variable exponent Sobolev space $W^{1,p(\cdot )}(\Omega )$ is defined by%
\begin{equation*}
W^{1,p(\cdot )}(\Omega )=\left\{ u\in L^{p(\cdot )}( \Omega )
;\ \nabla u\in \lbrack L^{p(\cdot )}( \Omega )
]^{N}\right\}
\end{equation*}%
and it is equipped with the norm
\begin{equation*}
\left\Vert u\right\Vert _{W^{1,p(\cdot )}(\Omega )}=\left\vert u\right\vert
_{L^{p(\cdot )}(\Omega )}+\left\vert \nabla u\right\vert _{L^{p(\cdot
)}(\Omega )},\ \forall u\in W^{1,p(\cdot )}\left( \Omega \right) .
\end{equation*}

We denote by $W_{0}^{1,p(\cdot )}(\Omega )$ the closure of $C_{0}^{\infty
}\left( \Omega \right) $ in $W^{1,p(\cdot )}\left( \Omega \right) $.

The Lebesgue and Sobolev spaces with variable exponents coincide with the
usual Lebesgue and Sobolev spaces provided that $p$ is constant. According
to [\cite[pp. 8-9]{radrep}], these function spaces $L^{p(\cdot )}$ and $%
W^{1,p(\cdot )}$ have some non-usual properties, such as:

(i) Assuming that $1<p^-\leq p^+<\infty$ and $p:\overline\Omega\rightarrow
[1,\infty)$ is a smooth function, then the following co-area formula
\begin{equation*}
\int_\Omega |u(x)|^pdx=p\int_0^\infty t^{p-1}\,|\{x\in\Omega ;\
|u(x)|>t\}|\,dt
\end{equation*}
has no analogue in the framework of variable exponents.

(ii) Spaces $L^{p(\cdot )}$ do \textit{not} satisfy the \textit{mean
continuity property}. More exactly, if $p$ is nonconstant and continuous in
an open ball $B$, then there is some $u\in L^{p(\cdot )}(B)$ such that $%
u(x+h)\not\in L^{p(\cdot )}(B)$ for every $h\in
\mathbb{R}
{}^{N}$ with arbitrary small norm.

(iii) Function spaces with variable exponent are \textit{never} invariant
with respect to translations. The convolution is also limited. For instance,
the classical Young inequality
\begin{equation*}
|f\ast g|_{p(\cdot )}\leq C\,|f|_{p(\cdot )}\,\Vert g\Vert _{L^{1}}
\end{equation*}%
remains true if and only if $p$ is constant.

\smallskip Conditions ($\mathcal{A}_{1}$)-(i) and (ii) imply that%
\begin{equation}
\mathscr{A}(x,\xi )\leq \mathbf{A}(x,\xi )\cdot \xi\ \text{ for all }(x,\xi
)\in
\mathbb{R}
^{N}\times
\mathbb{R}
^{N}.  \label{2.21}
\end{equation}

Furthermore, ($\mathcal{A}_{1}$)-(ii) is weaker than the request that $%
\mathscr{A}$ is uniformly convex, that is, for any $\varepsilon \in (0,1)$,
there exists a constant $\delta (\varepsilon )\in (0,1)$ such that%
\begin{equation*}
\mathscr{A}\left(x,\frac{\xi +\eta }{2}\right)\leq (1-\delta (\varepsilon ))\,\frac{%
\mathscr{A}(x,\xi )+\mathscr{A}(x,\eta )}{2}
\end{equation*}%
for all $x\in
\mathbb{R}
^{N}$ and $(\xi ,\eta) \in
\mathbb{R}
^{N}\times
\mathbb{R}
^{N}$ satisfy $\left\vert u-v\right\vert \geq \varepsilon \max \{\left\vert
u\right\vert ,\left\vert v\right\vert \}$.

By ($\mathcal{A}_{1}$)-(i) and (iii), we have%
\begin{equation*}
\mathscr{A}(x,\xi )=\int_{0}^{1}\frac{d}{dt}\mathscr{A}(x,t\xi
)dt=\int_{0}^{1}\frac{1}{t}\mathbf{A}(x,t\xi )\cdot t\xi dt\geq \left\{
\begin{array}{c}
c_{1}\left\vert \xi \right\vert ^{p(x)},\left\vert \xi \right\vert >1 \\
c_{1}\left\vert \xi \right\vert ^{q(x)},\left\vert \xi \right\vert \leq 1.
\end{array}%
\right.
\end{equation*}%
This estimate in combination with (\ref{5.1}) and (\ref{2.21}) yields
\begin{equation}
\left.
\begin{array}{c}
c_{1}\left\vert \xi \right\vert ^{p(x)},\left\vert \xi \right\vert >1 \\
c_{1}\left\vert \xi \right\vert ^{q(x)},\left\vert \xi \right\vert \leq 1%
\end{array}%
\right\} \leq \mathscr{A}(x,\xi )\leq \mathbf{A}(x,\xi )\cdot \xi \text{ }%
\leq \left\{
\begin{array}{c}
c_{2}\left\vert \xi \right\vert ^{p(x)},\left\vert \xi \right\vert >1 \\
c_{2}\left\vert \xi \right\vert ^{q(x)},\left\vert \xi \right\vert \leq 1%
\end{array}%
\right. ,\quad\forall (x,\xi )\in
\mathbb{R}
^{N}\times
\mathbb{R}
^{N}.  \label{2.1}
\end{equation}

Denote
\begin{equation*}
L_{w}^{\vartheta (\cdot )}(\Omega )=\left\{ u\left\vert \int_{\Omega
}w(x)\left\vert u(x)\right\vert ^{\vartheta (x)}dx<\infty \right. \right\} ,
\end{equation*}%
with the norm%
\begin{equation*}
\left\vert u\right\vert _{L_{w}^{\vartheta (\cdot )}(\Omega )}=\inf \left\{
\lambda >0\left\vert \int_{\Omega }w(x)\left\vert \frac{u(x)}{\lambda }%
\right\vert ^{\vartheta (x)}dx\leq 1\right. \right\} .
\end{equation*}

If $\Omega =%
\mathbb{R}
^{N}$, we simply denote ($L_{w}^{\vartheta (\cdot )}(%
\mathbb{R}
^{N}),\left\vert \cdot \right\vert _{L_{w}^{\vartheta (\cdot )}(%
\mathbb{R}
^{N})}$) as ($L_{w}^{\vartheta (\cdot )},\left\vert \cdot \right\vert
_{L_{w}^{\vartheta (\cdot )}}$).

From now on, we denote by $B_{R}$ the ball in $%
\mathbb{R}
^{N}$ centered at the origin and of radius $R>0$.\newline

\textbf{Lemma 2.5.} (see [\cite{w5}, Lemma 2.2]) Assume that $\vartheta ^{-}>1$
and $\vartheta ^{+}<\infty $. Then $L_{w}^{\vartheta (\cdot )}$ is separable
uniformly convex Banach space.\newline

\textbf{Theorem 2.6.} (see [\cite{w12}] Interpolation Theorem) If $p(\cdot
)<\alpha (\cdot )<q(\cdot )$, then for any $u\in L^{\alpha (\cdot )}(\Omega
) $, there exists $\lambda =\lambda (\Omega ,\alpha ,p,q,u)\in \lbrack
\theta ^{-},\theta ^{+}]$, where $\theta (\cdot )=\frac{p(q-\alpha )}{\alpha
(q-p)}$, such that
\begin{equation*}
\left\vert u\right\vert _{L^{\alpha (\cdot )}(\Omega )}\leq 2\left\vert
u\right\vert _{L^{p(\cdot )}(\Omega )}^{\lambda }\cdot \left\vert
u\right\vert _{L^{q(\cdot )}(\Omega )}^{1-\lambda }.
\end{equation*}%
Moreover, if $\theta ^{-}<\theta ^{+}$, then $\lambda \in (\theta
^{-},\theta ^{+})$.\newline

\textbf{Proposition 2.7 } (see [\cite{j5}, \cite{e12}]) If $\Omega $ is a
bounded domain, we have

(i) $W^{1,p(\cdot )}(\Omega )$ and $W_{0}^{1,p(\cdot )}(\Omega )$ are
separable reflexive Banach spaces;

(ii) if $\theta \in C_{+}\left( \overline{\Omega }\right) $ and $\theta
(x)<p^{\ast }(x)$ for any $x\in \overline{\Omega },$ then the imbedding from
$W^{1,p(\cdot )}(\Omega )$ to $L^{\theta (\cdot )}\left( \Omega \right) $ is
compact and continuous;

(iii) there is a constant $C>0,$ such that
\begin{equation*}
\left\vert u\right\vert _{L^{p(\cdot )}(\Omega )}\leq C\left\vert \nabla
u\right\vert _{L^{p(\cdot )}(\Omega )},\ \forall u\in W_{0}^{1,p(\cdot
)}(\Omega ).
\end{equation*}

\section{Variable Exponent Orlicz-Sobolev Spaces Theory}

Let $\Omega \subset
\mathbb{R}
^{N}$ be an open domain.\newline

\textbf{Definition 3.1.} We define the following real valued linear space
\begin{equation*}
L^{p(\cdot )}(\Omega )+L^{q(\cdot )}(\Omega )=\left\{ u\mid u=v+w,v\in
L^{p(\cdot )}(\Omega ),w\in L^{q(\cdot )}(\Omega )\right\} ,
\end{equation*}%
which is endowed with the norm
\begin{equation}
\left\vert u\right\vert _{L^{p(\cdot )}(\Omega )+L^{q(\cdot )}(\Omega
)}=\inf \left\{ \left\vert v\right\vert _{L^{p(\cdot )}(\Omega )}+\left\vert
w\right\vert _{L^{q(\cdot )}(\Omega )}\mid v\in L^{p(\cdot )}(\Omega ),w\in
L^{q(\cdot )}(\Omega ),v+w=u\right\} .  \label{o1}
\end{equation}

If $\Omega =%
\mathbb{R}
^{N}$, we simply denote ($L^{p(\cdot )}(\Omega )+L^{q(\cdot )}(\Omega
),\left\vert \cdot \right\vert _{L^{p(\cdot )}(\Omega )+L^{q(\cdot )}(\Omega
)}$) as ($L^{p(\cdot )}+L^{q(\cdot )},\left\vert \cdot \right\vert
_{L^{p(\cdot )}+L^{q(\cdot )}}$).

We also define the linear space
\begin{equation*}
L^{p(\cdot )}(\Omega )\cap L^{q(\cdot )}(\Omega )=\left\{ u\mid u\in
L^{p(\cdot )}(\Omega )\text{ and }u\in L^{q(\cdot )}(\Omega )\right\} ,
\end{equation*}%
which is endowed with the norm
\begin{equation*}
\left\vert u\right\vert _{L^{p(\cdot )}(\Omega )\cap L^{q(\cdot )}(\Omega
)}=\max \left\{ \left\vert u\right\vert _{L^{p(\cdot )}(\Omega )},\left\vert
u\right\vert _{L^{q(\cdot )}(\Omega )}\right\} .
\end{equation*}

Throughout this paper, we denote
\begin{equation*}
\Lambda _{u}=\left\{ x\in \Omega \mid\ \left\vert u(x)\right\vert >1\right\}
\quad \mbox{and}\quad \Lambda _{u}^{c}=\left\{ x\in \Omega \mid\ \left\vert
u(x)\right\vert \leq 1\right\} .
\end{equation*}

\textbf{Proposition 3.2. }(see [\cite{w12}, Proposition 3.2]) Assume ($%
\mathcal{A}_{1}$)-(iv). Let $\Omega \subset
\mathbb{R}
^{N}$ and $u\in L^{p(\cdot )}(\Omega )+L^{q(\cdot )}(\Omega )$. Then the
following properties hold:

(i) if $\Omega ^{\prime }\subset \Omega $ \ is such that $\left\vert \Omega
^{\prime }\right\vert <+\infty $, then $u\in L^{p(\cdot )}(\Omega ^{\prime
}) $;

(ii) if $\Omega ^{\prime }\subset \Omega $ \ is such that $u\in L^{\infty
}(\Omega ^{\prime })$, then $u\in L^{q(\cdot )}(\Omega ^{\prime })$;

(iii) $\left\vert \Lambda _{u}\right\vert <+\infty $;

(iv) $u\in L^{p(\cdot )}(\Lambda _{u})\cap L^{q(\cdot )}(\Lambda _{u}^{c})$;

(v) the infimum in (\ref{o1}) is attained;

(vi) if $B\subset \Omega $, then $\left\vert u\right\vert _{L^{p(\cdot
)}(\Omega )+L^{q(\cdot )}(\Omega )}\leq \left\vert u\right\vert _{L^{p(\cdot
)}(B)+L^{q(\cdot )}(B)}+\left\vert u\right\vert _{L^{p(\cdot )}(\Omega
/B)+L^{q(\cdot )}(\Omega /B)}$;

(vii) we have
\begin{equation*}
\begin{array}{ll}
& \displaystyle\max \left\{ \frac{1}{1+2\left\vert \Lambda _{u}\right\vert ^{%
\frac{1}{p(\xi )}-\frac{1}{q(\xi )}}}\left\vert u\right\vert _{L^{p(\cdot
)}(\Lambda _{u})},c\min \{\left\vert u\right\vert _{L^{q(\cdot )}(\Lambda
_{u}^{c})},\left\vert u\right\vert _{L^{q(\cdot )}(\Lambda _{u}^{c})}^{\frac{%
q(\xi )}{p(\xi )}}\}\right\} \leq \\
& \displaystyle\left\vert u\right\vert _{L^{p(\cdot )}(\Omega )+L^{q(\cdot
)}(\Omega )}\leq \left\vert u\right\vert _{L^{p(\cdot )}(\Lambda
_{u})}+\left\vert u\right\vert _{L^{q(\cdot )}(\Lambda _{u}^{c})},%
\end{array}%
\end{equation*}%
where $\xi \in
\mathbb{R}
^{N}$ and $c$ is a small positive constant.\newline

\textbf{Proposition 3.3.} (see [\cite{w11}, Theorem 2]) If ($X,\left\Vert
\cdot \right\Vert $) is a Banach space, then the following two statements
are equivalent:

(i) ($X,\left\Vert \cdot \right\Vert $) is reflexive;

(ii) any bounded sequence of ($X,\left\Vert \cdot \right\Vert $) has a weak
convergent subsequence.\newline

\textbf{Proposition 3.4.} (see [\cite{w12}, Proposition 3.8]) Assume that hypothesis ($%
\mathcal{A}_{1}$)-(iv) is fulfilled. Then $(L^{p^{\prime }(\cdot )}(\Omega )\cap
L^{q^{\prime }(\cdot )}(\Omega ))^{\prime }=L^{p(\cdot )}(\Omega
)+L^{q(\cdot )}(\Omega )$.\newline

\textbf{Proposition 3.5.} (see [\cite{w12}, Proposition 3.9]) Assume that hypothesis (%
$\mathcal{A}_{1}$)-(iv) is fulfilled. Then ($L^{p(\cdot
)}(\Omega )+L^{q(\cdot )}(\Omega ),\left\vert \cdot \right\vert _{L^{p(\cdot
)}(\Omega )+L^{q(\cdot )}(\Omega )}$) is a reflexive Banach space.\newline

Define $X(\Omega )\mathcal{=}\left\{ u\in L^{\alpha (\cdot )}(\Omega )\mid
\nabla u\in (L^{p(\cdot )}(\Omega )+L^{q(\cdot )}(\Omega ))^{N}\right\} $
with the following norm
\begin{equation*}
\left\Vert u\right\Vert _{\Omega }=\left\vert u\right\vert _{L^{\alpha
(\cdot )}(\Omega )}+\left\vert \nabla u\right\vert _{L^{p(\cdot )}(\Omega
)+L^{q(\cdot )}(\Omega )}.
\end{equation*}

If $\Omega =%
\mathbb{R}
^{N}$, we simply denote ($X(\Omega ),\left\Vert u\right\Vert _{\Omega }$) as
($X,\left\Vert u\right\Vert $).\newline

\textbf{Proposition 3.6.} (see [\cite{w12}, Proposition 3.10]) Assume ($%
\mathcal{A}_{1}$)-(iv). Then ($X(\Omega ),\left\Vert u\right\Vert _{\Omega }$%
) is a Banach space.\newline

\textbf{Proposition 3.7.} (see [\cite{w12}, Proposition 3.11]) Assume ($%
\mathcal{A}_{1}$)-(iv). Then ($X(\Omega ),\left\Vert u\right\Vert _{\Omega }$%
) is reflexive.\newline

\textbf{Theorem 3.8.} (see [\cite{w12}, Theorem 3.12]) Assume ($\mathcal{A}%
_{1}$)-(iv), $1<<p^{\ast }(\cdot )\frac{q^{\prime }(\cdot )}{p^{\prime
}(\cdot )}$, $\alpha $ satisfies $1<<\alpha (\cdot )<<p^{\ast }(\cdot )\frac{%
N-1}{N}$ and $1<<\alpha (\cdot )\leq p^{\ast }(\cdot )\frac{q^{\prime
}(\cdot )}{p^{\prime }(\cdot )}$. Then the space $X(\Omega )$ is
continuously embedded into $L^{p^{\ast }(\cdot )}(\Omega )$.\newline

\textbf{Corollary 3.9.} (see [\cite{w12}, Corollary 3.13]) Assume conditions
of Theorem 3.8. We have the following properties:

(i) for any $u\in X(\Omega )$, $\psi _{n}u\rightarrow u$ in $X(\Omega )$;

(ii) for any $u\in X$, we have $u_{\varepsilon }=u\ast \mathbf{j}%
_{\varepsilon }\rightarrow u$ in $X$ (where $\mathbf{j}_{\varepsilon
}(x)=\varepsilon ^{-N}\mathbf{j}(\frac{x}{\varepsilon })$ and $\mathbf{j}:%
\mathbb{R}
^{N}\rightarrow
\mathbb{R}
^{+}$ is in $C_{c}^{\infty }(%
\mathbb{R}
^{N})$, a function inducing a probability measure);

(iii) for any $u\in X$, there exists a sequence \{$u_{n}$\}$\subset
C_{c}^{\infty }(%
\mathbb{R}
^{N})$ such that $u_{n}\rightarrow u$ in $X$.\newline

\textbf{Theorem 3.10.} Assume conditions of Theorem 3.8.

(i) For any $\alpha \leq s\leq p^{\ast }$, the space $X(\Omega )$ is
continuously embedded into $L^{s(\cdot )}(\Omega )$.

(ii) For any bounded subset $\Omega \subset
\mathbb{R}
^{N}$, there is a compact embedding $X(\Omega )\mathcal{\hookrightarrow }%
L^{s(\cdot )}(\Omega )$ for any $1\leq s\leq p^{\ast }$.

(iii) We also assume that $\vartheta (\cdot )\in C(%
\mathbb{R}
^{N})$ is Lipschitz continuous, $w\in L^{r(\cdot )}$ and
\begin{equation}
\alpha (\cdot )\leq r^{\prime }(\cdot )\vartheta (\cdot )\leq p^{\ast
}(\cdot )\text{ in }%
\mathbb{R}
^{N}.  \label{w5}
\end{equation}%
Then there is a compact embedding $X\mathcal{\hookrightarrow }%
L_{w}^{\vartheta (\cdot )}$.

\smallskip \textbf{Proof}. The proofs of (i) and (ii) are trivial from
Proposition 2.7. We only need to prove (iii).

Since $X$ is  embedded into $L^{r^{\prime }(\cdot )\vartheta (\cdot )}(%
\mathbb{R}
^{N})$ for $\alpha (\cdot )\leq r^{\prime }(\cdot )\vartheta (\cdot
)<<p^{\ast }(\cdot )$, we may assume that $u_{n}\rightharpoonup u$ in $X$.
Then \{$\left\Vert u_{n}\right\Vert $\} is bounded and the continuous
embedding $X\hookrightarrow $ $L^{r^{\prime }(\cdot )\vartheta (\cdot )}(%
\mathbb{R}
^{N})$ guarantees the boundedness of \{$|u_{n}|_{L^{r^{\prime }(\cdot
)\vartheta (\cdot )}}$\}. So, there is a positive constant $M$ such that%
\begin{equation*}
\sup \{||u_{n}|^{\vartheta (x)}|_{L^{r^{\prime }(\cdot )}},||u|^{\vartheta
(x)}|_{L^{r^{\prime }(\cdot )}}\}\leq M.
\end{equation*}

Set $B_{k}=\{x\in
\mathbb{R}
^{N}\mid |x|<k\}$. If $w\in L^{r(\cdot )}(%
\mathbb{R}
^{N})$ then
\begin{equation*}
|w|_{L^{r(\cdot )}(%
\mathbb{R}
^{N}\backslash B_{k})}\rightarrow 0,\ \mbox{as}\text{ }k\rightarrow \infty \text{.}
\end{equation*}

For any $\varepsilon >0$, we can find large enough $k_{1}>0$  such
that
\begin{equation*}
|w|_{L^{r(\cdot )}(%
\mathbb{R}
^{N}\backslash B_{k})}\leq \frac{\varepsilon }{8M}\text{, }\quad\mbox{for all}\ k\geq k_{1}\text{.%
}
\end{equation*}

From (ii) of this theorem, there is a compact embedding $X(B_{k_{1}})%
\hookrightarrow $ $L^{r^{\prime }(\cdot )\vartheta (\cdot )}(B_{k_{1}})$, so
$u_{n}\rightharpoonup u$ implies%
\begin{equation*}
\int_{B_{k_{1}}}|w(x)||u_{n}-u|^{\vartheta (x)}dx\leq |w(x)|_{L^{r(\cdot
)}(B_{k_{1}})}||u_{n}-u|^{\vartheta (x)}|_{L^{r^{\prime }(\cdot
)}(B_{k_{1}})}\rightarrow 0\text{ as }n\rightarrow +\infty \text{.}
\end{equation*}

Thus, there exists $n_{1}>0$ such that for all $n\geq n_{1}$ we have
\begin{equation*}
\int_{B_{k_{1}}}|w(x)||u_{n}-u|^{\vartheta (x)}dx<\frac{\varepsilon }{2}%
\text{.}
\end{equation*}

Therefore
\begin{eqnarray*}
\int_{%
\mathbb{R}
^{N}}|w(x)||u_{n}-u|^{\vartheta (x)}dx &\leq
&\int_{B_{k_{1}}}|w(x)||u_{n}-u|^{\vartheta (x)}dx \\
&&+\int_{%
\mathbb{R}
^{N}\backslash B_{k_{1}}}|w(x)||u_{n}-u|^{\vartheta (x)}dx \\
&\leq &\frac{\varepsilon }{2}+2|w|_{L^{r(\cdot )}(%
\mathbb{R}
^{N}\backslash B_{k})}||u_{n}-u|^{\vartheta (x)}|_{L^{r^{\prime }(\cdot )}(%
\mathbb{R}
^{N}\backslash B_{k})} \\
&\leq &\frac{\varepsilon }{2}+2|w|_{L^{r(\cdot )}(%
\mathbb{R}
^{N}\backslash B_{k})}\{||u_{n}|^{\vartheta (x)}|_{L^{r^{\prime }(\cdot )}(%
\mathbb{R}
^{N}\backslash B_{k})}+||u|^{\vartheta (x)}|_{L^{r^{\prime }(\cdot )}(%
\mathbb{R}
^{N}\backslash B_{k})}\} \\
&\leq &\frac{\varepsilon }{2}+\frac{\varepsilon }{2}=\varepsilon .
\end{eqnarray*}

We conclude that $u_{n}\rightarrow u$ in $L_{w}^{\vartheta (\cdot )}$. This
completes the proof. $\square $

\section{Properties of Functionals and Operators}

By (vii) of Proposition 3.2, we deduce that $\mathscr{A}(x,\nabla u)$ is
integrable on $%
\mathbb{R}
^{N}$ for all $u\in X$. Thus, $\int_{%
\mathbb{R}
^{N}}\mathscr{A}(x,\nabla u)dx$ is well defined. For $u\in X$, it follows by
(\ref{2.1}) that%
\begin{eqnarray}
&&\int_{%
\mathbb{R}
^{N}}\mathbf{A}(x,\nabla u)\cdot \nabla udx+\int_{%
\mathbb{R}
^{N}}\left\vert u\right\vert ^{\alpha (x)}dx  \notag \\
&\geq &c_{1}\left( \int_{%
\mathbb{R}
^{N}\cap \Lambda _{\nabla u}}\left\vert \nabla u\right\vert ^{p(x)}dx+\int_{%
\mathbb{R}
^{N}\cap \Lambda _{\nabla u}^{c}}\left\vert \nabla u\right\vert
^{q(x)}dx+\int_{%
\mathbb{R}
^{N}}\left\vert u\right\vert ^{\alpha (x)}dx\right) ,  \label{2.5}
\end{eqnarray}%
and%
\begin{eqnarray*}
&&\int_{%
\mathbb{R}
^{N}}\mathbf{A}(x,\nabla u)\cdot \nabla udx+\int_{%
\mathbb{R}
^{N}}\left\vert u\right\vert ^{\alpha (x)}dx \\
&\leq &c_{2}\left( \int_{%
\mathbb{R}
^{N}\cap \Lambda _{\nabla u}}\left\vert \nabla u\right\vert ^{p(x)}dx+\int_{%
\mathbb{R}
^{N}\cap \Lambda _{\nabla u}^{c}}\left\vert \nabla u\right\vert
^{q(x)}dx+\int_{%
\mathbb{R}
^{N}}\left\vert u\right\vert ^{\alpha (x)}dx\right)
\end{eqnarray*}%
where $c_{1}$ and $c_{2}$ are positive constants.

Similarly, using (\ref{2.1}), we get for all $u\in X$
\begin{eqnarray}
&&\int_{%
\mathbb{R}
^{N}}\mathscr{A}(x,\nabla u)dx+\int_{%
\mathbb{R}
^{N}}\frac{1}{\alpha (x)}\left\vert u\right\vert ^{\alpha (x)}dx  \notag \\
&\geq &c_{1}\left( \int_{%
\mathbb{R}
^{N}\cap \Lambda _{\nabla u}}\left\vert \nabla u\right\vert ^{p(x)}dx+\int_{%
\mathbb{R}
^{N}\cap \Lambda _{\nabla u}^{c}}\left\vert \nabla u\right\vert
^{q(x)}dx+\int_{%
\mathbb{R}
^{N}}\frac{1}{\alpha (x)}\left\vert u\right\vert ^{\alpha (x)}dx\right) ,
\label{2.52}
\end{eqnarray}%
and
\begin{eqnarray*}
&&\int_{%
\mathbb{R}
^{N}}\mathscr{A}(x,\nabla u)dx+\int_{%
\mathbb{R}
^{N}}\frac{1}{\alpha (x)}\left\vert u\right\vert ^{\alpha (x)}dx \\
&\leq &c_{2}\left( \int_{%
\mathbb{R}
^{N}\cap \Lambda _{\nabla u}}\left\vert \nabla u\right\vert ^{p(x)}dx+\int_{%
\mathbb{R}
^{N}\cap \Lambda _{\nabla u}^{c}}\left\vert \nabla u\right\vert
^{q(x)}dx+\int_{%
\mathbb{R}
^{N}}\frac{1}{\alpha (x)}\left\vert u\right\vert ^{\alpha (x)}dx\right) .
\end{eqnarray*}

From ($\mathcal{H}_{f}^{1}$), we have $\left\vert g(x,u)\right\vert \leq
\left\vert u\right\vert ^{\alpha (x)-1}+C\left\vert u\right\vert ^{\gamma
(x)-1}$. Notice that $\alpha (\cdot )\leq \gamma (\cdot )\leq \frac{p^{\ast
}(\cdot )}{r^{\prime }(\cdot )}$. Combining ($\mathcal{H}_{w}$) and ($%
\mathcal{H}_{a}$), it is easy to check that $f(x,u)v$ and $F(x,u)$ are
integrable on $%
\mathbb{R}
^{N}$ for all $u,v\in X$.

We say that $u\in X$ is a solution of problem $(\mathcal{E})$ if%
\begin{equation*}
\int_{%
\mathbb{R}
^{N}}\mathbf{A}(x,\nabla u)\cdot \nabla vdx+\int_{%
\mathbb{R}
^{N}}\left\vert u\right\vert ^{\alpha (x)-2}uvdx=\int_{%
\mathbb{R}
^{N}}f(x,u)vdx,\quad \forall v\in X.
\end{equation*}

It follows that solutions of $(\mathcal{E})$ correspond to the critical
points of the Euler-Lagrange energy functional $\Phi :X\rightarrow
\mathbb{R}
$, defined by%
\begin{equation*}
\Phi =\int_{%
\mathbb{R}
^{N}}\mathscr{A}(x,\nabla u)dx+\int_{%
\mathbb{R}
^{N}}\frac{1}{\alpha (x)}\left\vert u\right\vert ^{\alpha (x)}dx-\int_{%
\mathbb{R}
^{N}}F(x,u)dx,
\end{equation*}%
where $F(x,u)=\int_{0}^{u}f(x,s)ds$.

Define functionals $\Phi _{\mathscr{A}},\Phi _{\alpha },\Phi
_{f}:X\rightarrow
\mathbb{R}
$ by
\begin{equation*}
\Phi _{\mathscr{A}}(u)=\int_{%
\mathbb{R}
^{N}}\mathscr{A}(x,\nabla u)dx\text{, }\ \Phi _{\alpha }(u)=\int_{%
\mathbb{R}
^{N}}\frac{1}{\alpha (x)}\left\vert u\right\vert ^{\alpha (x)}dx\text{, }\
\Phi _{f}(u)=\int_{%
\mathbb{R}
^{N}}F(x,u)dx.
\end{equation*}

\textbf{Lemma 4.1.} (see [\cite{w12}, Lemma 4.1]) Assume the structure
conditions ($\mathcal{A}_{1}$). Then the functional $\Phi _{\mathscr{A}}\ $%
is convex, of class $C^{1}$, and sequentially weakly lower semicontinuous in
$X$. Moreover, $\Phi _{\mathscr{A}}^{\prime }:X\rightarrow X^{\ast }$ is
bounded.\newline

\textbf{Lemma 4.2.} (see [\cite{w12}, Lemma 4.2]) Assume the structure
conditions ($\mathcal{A}_{1}$)-(iv). Then the functional $\Phi _{\alpha }\ $%
is convex, of class $C^{1}$ and sequentially weakly lower semicontinuous.
Moreover, if $u_{n},u\in X$ and $u_{n}\rightharpoonup u\ $in $X$, then $\Phi
_{\alpha }^{\prime }(u_{n})\overset{\ast }{\rightharpoonup }\Phi _{\alpha
}^{\prime }(u)$ in $X^{\ast }$.\newline

\textbf{Lemma 4.3.} (see [\cite{w12}, Lemma 4.3]) Assume ($\mathcal{A}_{1}$%
), ($\mathcal{H}_{f}^{1}$), ($\mathcal{H}_{w}$) and ($\mathcal{H}_{a}$).
Then $\Phi _{f}$ is of class $C^{1}$ and sequentially weakly-strongly
continuous, that is, if $u_{n}\rightharpoonup u$ in $X$ then $\Phi
_{f}(u_{n})\rightarrow \Phi _{f}(u)$ and $\Phi _{f}^{\prime
}(u_{n})\rightarrow \Phi _{f}^{\prime }(u)$ in $X^{\ast }$.

\textbf{Proof}. Since $X$ is embedded into $L^{\gamma (\cdot )}(%
\mathbb{R}
^{N})$ for $\alpha (\cdot )\leq \gamma (\cdot )\leq \frac{p^{\ast }(\cdot )}{%
r^{\prime }(\cdot )}$, we deduce that $F(x,u)$ is integrable on $%
\mathbb{R}
^{N}$, hence $\Phi _{f}(u)$ is well defined.

Now, let us prove that is $\Phi _{f}$ weakly-strongly continuous. Assume
that $u_{n}\rightharpoonup u$ in $X$, then \{$u_{n}$\} is bounded in $X$,
hence \{$\left\vert u_{n}\right\vert _{L^{r^{\prime }(\cdot )\alpha
(\cdot )}}$\} and \{$\left\vert u_{n}\right\vert _{L^{r^{\prime }(\cdot
)\gamma (\cdot )}}$\} are bounded.

Since%
\begin{equation*}
|G(x,t)|\leq \frac{1}{\alpha (x)}|t|^{\alpha (x)}+\frac{c}{\theta (x)}%
|t|^{\gamma (x)}
\end{equation*}%
we have
\begin{equation*}
w(x)|G(x,u_{n})-G(x,u)|\leq w(x)\{|u|^{\alpha (x)}+|u_{n}|^{\alpha
(x)}+c|u_{n}|^{\gamma (x)}+c|u|^{\gamma (x)}\}.
\end{equation*}

Therefore \{$w(x)|G(x,u_{n})-G(x,u)|$\} is uniformly integrable in $%
\mathbb{R}
^{N}$.

By Theorem 3.10, we have $u_{n}\rightarrow u$ a.e. in $%
\mathbb{R}
^{N}$.

Thus, by Vitali's theorem, we have%
\begin{equation*}
\underset{n\rightarrow \infty }{\lim }\int_{%
\mathbb{R}
^{N}}w(x)|G(x,u_{n})-G(x,u)|dx=\int_{%
\mathbb{R}
^{N}}\underset{n\rightarrow \infty }{\lim }w(x)|G(x,u_{n})-G(x,u)|dx=0.
\end{equation*}

Similarly, we have
\begin{equation*}
\underset{n\rightarrow \infty }{\lim }\int_{%
\mathbb{R}
^{N}}\frac{a(x)}{\delta (x)}|\left\vert u_{n}\right\vert ^{\delta
(x)}-\left\vert u\right\vert ^{\delta (x)}|dx=\int_{%
\mathbb{R}
^{N}}\underset{n\rightarrow \infty }{\lim }\frac{a(x)}{\delta (x)}%
|\left\vert u_{n}\right\vert ^{\delta (x)}-\left\vert u\right\vert ^{\delta
(x)}|dx=0.
\end{equation*}

We conclude that $\Phi _{f}(u_{n})\rightarrow \Phi _{f}(u)$.

In a similar way, we can obtain the weakly--strongly continuity of $\Phi
_{f}^{\prime }$. $\square $\newline

\textbf{Lemma 4.4.} Assume ($\mathcal{A}_{1}$), ($\mathcal{H}_{f}^{1}$), ($%
\mathcal{H}_{w}$) and ($\mathcal{H}_{a}$). Then the functional $\Phi $ is of
class $C^{1}$ and sequentially weakly lower semicontinuous in $X$, that is,
if $u_{n}\rightharpoonup u_{0}$ in $X$, then
\begin{equation*}
\Phi (u_{0})\leq \underset{n\rightarrow \infty }{\lim \inf }\Phi (u_{n}).
\end{equation*}

\textbf{Proof.} According to Lemmas 4.1-4.3, we deduce the $C^{1}$
continuity of $\Phi $. Next, we will prove that $\Phi $ is the sequentially
weakly lower semicontinuous in $X$.

By Lemma 4.3, $\Phi _{f}(u)$ is weakly continuous. Obviously%
\begin{eqnarray*}
\underset{n\rightarrow \infty }{\lim \inf }\Phi (u_{n}) &\geq &\underset{%
n\rightarrow \infty }{\lim \inf }(\Phi _{\mathscr{A}}(u_{n})+\Phi
_{a}(u_{n}))-\underset{n\rightarrow \infty }{\lim \sup }\Phi _{f}(u_{n}) \\
&\geq &\Phi _{\mathscr{A}}(u_{0})+\Phi _{\alpha }(u_{0})-\Phi _{f}(u_{0}) \\
&=&\Phi (u_{0}).
\end{eqnarray*}

Thus, $\Phi $ is sequentially weakly lower semicontinuous in $X$. $\square $%
\newline

\textbf{Lemma 4.5 }(see [\cite{w12}, Lemma 4.6]) Suppose that $\mathscr{A}$
satisfies ($\mathcal{A}_{1}$) and ($\mathcal{A}_{2}$) (namely $\mathscr{A}%
(x,\cdot ):%
\mathbb{R}
^{N}\rightarrow
\mathbb{R}
$ is a uniformly convex function), that is, for any $\varepsilon \in (0,1)$
there exists $\delta (\varepsilon )\in (0,1)$ such that  $\mathscr{A}\left(x,\frac{u+v}{2}\right)\leq \frac{(1-\delta
(\varepsilon ))}{2}(\mathscr{A}(x,u)+\mathscr{A}(x,v))$ for all $x\in {\mathbb R}^N$ and all $(u,v)\in {\mathbb R}^N$ with $\left\vert
u-v\right\vert \leq \varepsilon \max \{\left\vert u\right\vert ,\left\vert
v\right\vert \}$.
Then we have

(i) $\Phi _{\mathscr{A}}$($\cdot $) $:X\rightarrow
\mathbb{R}
$ is uniformly convex, that is, for any $\varepsilon \in (0,1)$ there exists
$\delta (\varepsilon )\in (0,1)$ such that for all $u,v\in X$
\begin{equation*}
\Phi _{\mathscr{A}}\left(\frac{u-v}{2}\right)\leq \varepsilon \frac{\Phi _{\mathscr{A}%
}(u)+\Phi _{\mathscr{A}}(v)}{2}\text{ or }\Phi _{\mathscr{A}}\left(\frac{u+v}{2}%
\right)\leq (1-\delta (\varepsilon ))\frac{\Phi _{\mathscr{A}}(u)+\Phi _{%
\mathscr{A}}(v)}{2};
\end{equation*}

(ii) if $u_{n}\rightharpoonup u$ in $X$ and $\underset{n\rightarrow \infty }{%
\overline{\lim }}$ $(\Phi _{\mathscr{A}}^{\prime }(u_{n})-\Phi _{\mathscr{A}%
}^{\prime }(u),u_{n}-u)\leq 0$, then $\Phi _{\mathscr{A}}(u_{n}-u)%
\rightarrow 0$ and $\left\vert \nabla u_{n}-\nabla u\right\vert _{L^{p(\cdot
)}+L^{q(\cdot )}}\rightarrow 0$.\newline

Define $\rho (\cdot ):X\rightarrow
\mathbb{R}
$ as%
\begin{equation*}
\rho (u)=\int_{%
\mathbb{R}
^{N}}\mathscr{A}(x,\nabla u)dx+\int_{%
\mathbb{R}
^{N}}\frac{1}{\alpha (x)}\left\vert u\right\vert ^{\alpha (x)}dx,
\end{equation*}%
and we denote the derivative operator by $L$, that is, $L=\rho ^{\prime
}:X\rightarrow X^{\ast }$ with
\begin{equation*}
(L(u),v)=\int_{%
\mathbb{R}
^{N}}\mathbf{A}(x,\nabla u)\nabla vdx+\int_{%
\mathbb{R}
^{N}}\left\vert u\right\vert ^{\alpha (x)-2}uvdx\quad \forall u,v\in X.
\end{equation*}

\textbf{Lemma 4.6. }(see [\cite{w12}, Lemma 4.7]) Under the structure
conditions ($\mathcal{A}_{1}$), we have the following properties.

(i) $L:X\rightarrow X^{\ast }$ is a continuous, bounded and strictly monotone
operator.

If ($\mathcal{A}_{2}$) is also satisfied, we have

(ii) $L$ is a mapping of type $(S_{+})$, that is, if $u_{n}\rightharpoonup u$
in $X$ and $\underset{n\rightarrow \infty }{\overline{\lim }}$ $%
(L(u_{n})-L(u),u_{n}-u)\leq 0,$ then $u_{n}\rightarrow u$ in $X$;

(iii) $L:X\rightarrow X^{\ast }$ is a homeomorphism.\newline

\textbf{Lemma 4.7. }We assume the structure conditions ($\mathcal{A}_{1}$)-($%
\mathcal{A}_{2}$), ($\mathcal{H}_{f}^{1}$)-($\mathcal{H}_{f}^{2}$), ($%
\mathcal{H}_{w}$), ($\mathcal{H}_{a}$), $1<<\alpha (\cdot )<<p^{\ast }(\cdot
)\frac{q^{\prime }(\cdot )}{p^{\prime }(\cdot )}$ and $\alpha \leq p$. Then $%
\Phi $ satisfies the (PS) condition, that is, if $\{u_{n}\}\subset X$
satisfies $\Phi (u_{n})\rightarrow c$ and $\left\Vert \Phi ^{\prime
}(u_{n})\right\Vert _{X^{\ast }}\rightarrow 0$, then $\{u_{n}\}$ has a
convergent subsequence.

\smallskip \textbf{Proof.} Assume that $\{u_{n}\}$ is bounded. Then up to a
subsequence, $u_{n}\rightharpoonup u_{0}$. By Lemma 4.3, again up to a
subsequence, we have $\Phi _{f}^{\prime }(u_{n})\rightarrow \Phi
_{f}^{\prime }(u_{0})$ in $X^{\ast }$. By Lemma 4.6, $L^{-1}$ is continuous
from $X^{\ast }$ to $X$, hence $u_{n}\rightarrow L^{-1}\circ \Phi
_{f}^{\prime }(u_{0})$ in $X$.

We only need to prove that $\{u_{n}\}$ is bounded in $X$.

We argue by contradiction. Suppose not, then there exist $c\in $ $%
\mathbb{R}
$ and \{$u_{n}$\} $\subset X$ satisfying:%
\begin{equation*}
\Phi (u_{n})\rightarrow c,\ \left\Vert \Phi ^{\prime }(u_{n})\right\Vert
_{X^{\ast }}\rightarrow 0,\ \left\Vert u_{n}\right\Vert \rightarrow +\infty .
\end{equation*}

Since $(\Phi ^{\prime }(u_{n}),\frac{1}{\theta }u_{n})\rightarrow 0$, we may
assume that%
\begin{eqnarray*}
c+\left\Vert u_{n}\right\Vert &\geq &\Phi (u_{n})-(\Phi ^{\prime }(u_{n}),%
\frac{1}{\theta }u_{n}) \\
&=&\int_{%
\mathbb{R}
^{N}}\mathscr{A}(x,\nabla u_{n})dx+\int_{%
\mathbb{R}
^{N}}\frac{1}{\alpha (x)}\left\vert u_{n}\right\vert ^{\alpha (x)}dx-\int_{%
\mathbb{R}
^{N}}F(x,u_{n})dx \\
&&-\left\{\int_{%
\mathbb{R}
^{N}}\frac{1}{\theta }\mathbf{A}(x,\nabla u_{n})\nabla u_{n}dx+\int_{%
\mathbb{R}
^{N}}\frac{1}{\theta }\left\vert u_{n}\right\vert ^{\alpha (x)}dx-\int_{%
\mathbb{R}
^{N}}\frac{1}{\theta }f(x,u_{n})u_{n}dx\right\} \\
&\geq &\int_{%
\mathbb{R}
^{N}}\left(1-\frac{s(x)}{\theta }\right)\mathscr{A}(x,\nabla u_{n})dx \\
&&+\int_{%
\mathbb{R}
^{N}}\left(\frac{1}{\alpha (x)}-\frac{1}{\theta }\right)\left\vert u_{n}\right\vert
^{\alpha (x)}dx+\int_{%
\mathbb{R}
^{N}}\left\{\frac{1}{\theta }f(x,u_{n})u_{n}-F(x,u_{n})\right\}dx.
\end{eqnarray*}

It follows that
\begin{eqnarray*}
c+\left\Vert u_{n}\right\Vert &\geq &\Phi (u_{n})-(\Phi ^{\prime }(u_{n}),%
\frac{1}{\theta }u_{n}) \\
&\geq &c_{1}\int_{%
\mathbb{R}
^{N}}\mathscr{A}(x,\nabla u_{n})+\left\vert u_{n}\right\vert ^{\alpha
(x)}dx-\int_{%
\mathbb{R}
^{N}}\left(\frac{1}{\delta (x)}-\frac{1}{\theta }\right)\lambda a(x)\left\vert
u_{n}\right\vert ^{\delta (x)}dx \\
&\geq &c_{1}\left\Vert u_{n}\right\Vert ^{\alpha ^{-}}-C\left\Vert
u_{n}\right\Vert ^{\delta ^{+}}-C.
\end{eqnarray*}

Notice that $\alpha ^{-}>\delta ^{+}>1$. Thus, we obtain a contradiction. The proof of
Lemma 4.7 is complete. $\square $\newline

\section{Proof of Theorems}

In this section, we will give the proofs of Theorem 1.1-1.3.

\subsection{Proof of Theorem 1.1}

\textbf{Proof of Theorem 1.1 }Let us consider the following auxiliary
problem:
\begin{equation*}
-\mathrm{div}\,\mathbf{A}(x,\nabla u)+\left\vert u\right\vert ^{\alpha
(x)-2}u=f^{+}(x,u),\text{ \ }\eqno{(\mathcal{E}^{+})}
\end{equation*}%
where
\begin{equation*}
f^{+}(x,u)=\left\{
\begin{array}{c}
f(x,u),\quad\text{ if }f(x,u)\geq 0 \\
0,\quad\text{ if }f(x,u)<0.%
\end{array}%
\right.
\end{equation*}

The corresponding Euler-Lagrange functional is
\begin{equation*}
\Phi ^{+}(u)=\int_{%
\mathbb{R}
^{N}}\mathscr{A}(x,\nabla u)dx+\int_{%
\mathbb{R}
^{N}}\frac{1}{\alpha (x)}\left\vert u\right\vert ^{\alpha (x)}dx-\int_{%
\mathbb{R}
^{N}}F^{+}(x,u)dx,
\end{equation*}%
where $F^{+}(x,u)=\int_{0}^{u}f^{+}(x,t)dt$.

Similar to the proof of Lemma 4.7, we deduce that $\Phi ^{+}$ satisfies (PS)
condition.

Next, we prove that $\Phi ^{+}(u)$ satisfies the conditions of the mountain pass
lemma.

By assumption ($\mathcal{H}_{f}^{2}$), we have
\begin{equation*}
g(x,tu)tu\geq \theta G(x,tu)>0,\forall u\neq 0,\quad\forall t\neq 0,\ \forall x\in
\mathbb{R}
^{N},
\end{equation*}%
and%
\begin{equation*}
\frac{g(x,tu)u}{G(x,tu)}\geq \frac{\theta }{t}>0,\quad\forall u\neq 0,\ \forall
t>0,\forall x\in
\mathbb{R}
^{N}.
\end{equation*}

Integrating about $t$ from $1$ to $t$, we have%
\begin{equation}
G(x,tu)\geq \left\vert t\right\vert ^{\theta }G(x,u)\geq 0,\quad\forall u\in
\mathbb{R}
,\ \forall t\geq 1,\ \forall x\in
\mathbb{R}
^{N},  \label{w22}
\end{equation}%
\begin{equation*}
0\leq G(x,tu)\leq \left\vert t\right\vert ^{\theta }G(x,u),\quad\forall u\in
\mathbb{R}
,\ \forall t\in (0,1],\ \forall x\in
\mathbb{R}
^{N}.
\end{equation*}

Hence, if $0<\left\Vert u\right\Vert \leq 1$, then
\begin{eqnarray*}
\Phi ^{+}(u) &\geq &\frac{1}{s^{+}}\left\Vert u\right\Vert
^{s^{+}}-\left\Vert u\right\Vert ^{\theta }\int_{%
\mathbb{R}
^{N}}\mu w(x)G^{+}(x,\frac{u}{\left\Vert u\right\Vert })dx-\int_{%
\mathbb{R}
^{N}}\lambda \frac{a(x)}{\delta (x)}\left\vert u\right\vert ^{\delta (x)}dx
\\
&\geq &\frac{1}{s^{+}}\left\Vert u\right\Vert ^{s^{+}}-c\left\Vert
u\right\Vert ^{\theta }-\lambda \left\Vert u\right\Vert ^{\delta (\xi )},\
\text{ for some }\xi \in
\mathbb{R}
^{N}.
\end{eqnarray*}

Let $\epsilon >0$ be small enough and $\lambda >0$  small enough. Then
\begin{equation}
\Phi (u)\geq c>0\text{ for all }\left\Vert u\right\Vert =\epsilon .
\label{w21}
\end{equation}

For $0\leq u\in X\backslash \left\{ 0\right\} $ and $t>1$, we have%
\begin{eqnarray*}
\Phi ^{+}(tu) &=&\int_{%
\mathbb{R}
^{N}}\mathscr{A}(x,\nabla tu)dx+\int_{%
\mathbb{R}
^{N}}\frac{1}{\alpha (x)}\left\vert tu\right\vert ^{\alpha (x)}dx-\int_{%
\mathbb{R}
^{N}}F(x,tu)dx \\
&\leq &t^{s^{+}}\int_{%
\mathbb{R}
^{N}}\mathscr{A}(x,\nabla tu)dx+t^{\alpha ^{+}}\int_{%
\mathbb{R}
^{N}}\frac{1}{\alpha (x)}\left\vert tu\right\vert ^{\alpha (x)}dx-t^{\theta
}\int_{%
\mathbb{R}
^{N}}\mu w(x)G(x,u)dx.
\end{eqnarray*}

Since $\alpha ^{+}\leq s^{+}<\theta $, then $\Phi (tw)\rightarrow -\infty $ $%
(t\rightarrow +\infty )$. Obviously, $\Phi ^{+}\left( 0\right) =0$, then $%
\Phi ^{+}$ satisfies the conditions of the mountain pass lemma. So, $\Phi ^{+}$
admits at least one nontrivial critical point $u_{1}$ satisfies $\Phi
^{+}(u_{1})>0$.

Thus $(\mathcal{E}^{+})$ has a solution $u_{1}$, and it is easy to see that $%
u_{1}\geq 0$, so $u_{1}$ is a nonnegative solution $u_{1}$ of $(\mathcal{E})$
with $\Phi (u_{1})>0$.

Similarly, we can establish the existence of a nonpositive solution $u_{2}$
of $(\mathcal{E})$ with $\Phi (u_{2})>0$.

Define $h_{0}\in C_{0}(\overline{B(x_{0},\varepsilon _{0})})$ as
\begin{equation*}
h_{0}(x)=\left\{
\begin{array}{cc}
0, & \left\vert x-x_{0}\right\vert \geq \varepsilon _{0} \\
\varepsilon _{0}-\left\vert x-x_{0}\right\vert , & \left\vert
x-x_{0}\right\vert <\varepsilon _{0}.%
\end{array}%
\right.
\end{equation*}

Let $\varepsilon _{0}>0$ be small enough. By ($\mathcal{H}_{a}$), we have $%
\Phi ^{+}(th_{0})<0$ for small enough $t>0$. Combining (\ref{w21}) and Lemma
4.4, we deduce that $\Phi ^{+}$  attains its infimum in \{$u\in X\mid \left\Vert
u\right\Vert <\epsilon $\}. Therefore, $\Phi ^{+}$ admits at least one
nontrivial critical point $u_{3}$ satisfying $\Phi ^{+}(u_{3})<0$. It is
easy to see that $u_{3}\geq 0$, so $u_{3}$ is a nonnegative solution of $(%
\mathcal{E})$ with $\Phi (u_{3})<0$.

Similarly, we can establish the existence of a nonpositive solution $u_{4}$
of $(\mathcal{E})$ with $\Phi (u_{4})<0$. The proof is complete. $\square $

\subsection{Proof of Theorem 1.2}

In order to prove Theorem 1.2, we need to recall some preliminary results.
Since $X$ is a reflexive and separable Banach space (see [\cite{22}],
Section 17, Theorems 2-3), there exist sequences $\left\{ e_{j}\right\}
\subset X$ and $\left\{ e_{j}^{\ast }\right\} \subset X^{\ast }$ such that%
\begin{equation*}
X=\overline{\mathrm{span}}\,\{e_{j}\text{, }j=1,2,\cdots \}\text{, }\left.
{}\right. X^{\ast }=\overline{\mathrm{span}}^{w^{\ast }}\{e_{j}^{\ast }\text{%
, }j=1,2,\cdots \},
\end{equation*}%
and
\begin{equation*}
\langle e_{j}^{\ast },e_{j}\rangle=\left\{
\begin{array}{c}
1,\quad\mbox{if}\ i=j, \\
0,\quad\mbox{if}\ i\neq j.%
\end{array}%
\right.
\end{equation*}

For convenience, we write%
\begin{equation}
X_{j}=\mathrm{span}\,\{e_{j}\}\text{, }\ Y_{k}=\overset{k}{\underset{j=1}{%
\oplus }}X_{j}\text{, }\ Z_{k}=\overline{\overset{\infty }{\underset{j=k}{%
\oplus }}X_{j}}.\newline
\label{8.1}
\end{equation}

Let $\Theta (u)=\left\vert u\right\vert _{L_{a}^{\delta (\cdot
)}}+\left\vert u\right\vert _{L_{w}^{\alpha (\cdot )}}+\left\vert
u\right\vert _{L_{w}^{\gamma (\cdot )}}$. By Theorem 3.10, similar to the
proof of Lemma 4.3, we deduce that $\Theta :X\rightarrow
\mathbb{R}
$ is weakly-strongly continuous and $\Theta (0)=0$.\newline

\textbf{Lemma 5.1} (see [\cite{w5}, Lemma 5.1]) Assume that $\Theta
:X\rightarrow
\mathbb{R}
$ is weakly-strongly continuous and $\Theta (0)=0$, $\gamma _{0}>0$ is a
given number.\ Let
\begin{equation*}
\beta _{k}=\beta _{k}(\gamma _{0})=\sup \left\{ \Theta (u)\mid \left\Vert
u\right\Vert \leq \gamma _{0},u\in Z_{k}\right\} .
\end{equation*}%
Then $\beta _{k}\rightarrow 0$ as $k\rightarrow \infty $.\newline

To complete the proof of Theorem 1.2, we recall the following critical point
lemma (see, e.g., [\cite[Theorem 4.7]{48}]).\newline

\textbf{Lemma 5.2} Suppose that $\Phi \in C^{1}(X,R)$ is even and satisfies
the (PS) condition. Let $V^{+}$, $V^{-}\subset X$ be closed subspaces of $X$
with $\mathrm{codim}\,V^{+}+1=\mathrm{dim}\,V^{-}$, and suppose that the
following conditions are fulfilled:

($1^{0}$) $\Phi (0)=0$;

($2^{0}$) $\exists \tau >0,$ $\gamma _{0}>0$ such that $\forall u\in V^{+}:$
$\Vert u\Vert =\gamma _{0}\Rightarrow \Phi (u)\geq \tau $;

($3^{0}$) $\exists \rho >0\ $such that $\forall u\in V^{-}:$ $\Vert u\Vert
\geq \rho \Rightarrow \Phi (u)\leq 0.$

Consider the following set:
\begin{equation*}
\Gamma =\{h\in C^{0}(X,X)\mid h\text{ is odd, }h(u)=u\text{ if }u\in V^{-}%
\text{ and }\Vert u\Vert \geq \rho \}.
\end{equation*}%
Then

($a$) $\forall \delta _{0}>0$, $h\in \Gamma $, $S_{\delta _{0}}^{+}\cap
h(V^{-})\neq \varnothing $, where $S_{\delta _{0}}^{+}=\{u\in V^{+}\mid \Vert
u\Vert =\delta _{0}\};$

($b$) the number $\varpi :=\underset{h\in \Gamma }{\inf }\underset{\text{ }%
u\in V^{-}}{\sup }\Phi (h(u))\geq \tau >0$ is a critical value for $\Phi $.%
\newline

\textbf{Proof of Theorem 1.2} According to our assumptions, $\Phi $ is an
even functional and satisfies the (PS) compactness condition. Let $%
V_{k}^{+}=Z_{k}$, which is a closed linear subspace of $X$ and $%
V_{k}^{+}\oplus Y_{k-1}=X$.

Set $V_{k}^{-}=X_{k}$. We will prove that there are infinitely many pairs of
$V_{k}^{+}$ and $V_{k}^{-}$, such that $\varphi $ satisfies the conditions
of Lemma 5.2. We also show that the corresponding critical value $\varpi
_{k}:=\underset{h\in \Gamma }{\inf }\underset{\text{ }u\in V_{k}^{-}}{\sup }%
\Phi (h(u))$ tends to $+\infty $ as $k\rightarrow \infty $, which implies
that there are infinitely many pairs of solutions to the problem $(\mathcal{E%
})$.

For any $k=1,2,\cdots $, we prove that there exist $\rho _{k}>\tau _{k}>0$
and large enough $k$ such that
\begin{eqnarray*}
(A_{1})\qquad\text{ }b_{k} &:&=\inf \left\{ \Phi (u)\mid u\in
V_{k}^{+},\left\Vert u\right\Vert =\tau _{k}\right\} \rightarrow +\infty
\quad\mbox{as}\ k\rightarrow +\infty ; \\
(A_{2})\qquad\text{ }a_{k} &:&=\max \left\{ \Phi (u)\right\vert \text{ }u\in
V_{k}^{-},\left\Vert u\right\Vert =\rho _{k}\}\leq 0.
\end{eqnarray*}

We first show that ($A_{1}$) holds. Let $\sigma \in (0,1)$ be small enough.
By ($\mathcal{H}_{f}^{1}$), there exists $C(\sigma )>0$ such that
\begin{equation*}
G(x,u)\leq \sigma \left\vert u\right\vert ^{\alpha (x)}+C(\sigma )\left\vert
u\right\vert ^{\gamma (x)},\quad \forall x\in \mathbb{R}^{N},\ \forall u\in
\mathbb{R}.
\end{equation*}%
By computation, for any $u\in Z_{k}$ with $\left\Vert u\right\Vert =\tau
_{k}=(2C_{2}\frac{1}{c_{1}}\beta _{k}^{\delta ^{+}})^{1/(\alpha ^{-}-\gamma
^{+})}$, we have
\begin{eqnarray*}
\Phi (u) &=&\int_{%
\mathbb{R}
^{N}}\mathscr{A}(x,\nabla u)dx+\int_{%
\mathbb{R}
^{N}}\frac{1}{\alpha (x)}\left\vert u\right\vert ^{\alpha (x)}dx-\int_{%
\mathbb{R}^{N}}F(x,u)dx \\
&\geq &2c_{1}\left(\int_{%
\mathbb{R}
^{N}\cap \Lambda _{\nabla u}}\left\vert \nabla u\right\vert ^{p(x)}dx+\int_{%
\mathbb{R}
^{N}\cap \Lambda _{\nabla u}^{c}}\left\vert \nabla u\right\vert
^{q(x)}dx+\int_{%
\mathbb{R}
^{N}}\left\vert u\right\vert ^{\alpha (x)}dx\right) \\
&&-\int_{%
\mathbb{R}
^{N}}\lambda \frac{a(x)}{\delta (x)}\left\vert u\right\vert ^{\delta
(x)}dx-\sigma \int_{\mathbb{R}^{N}}\mu w(x)\left\vert u\right\vert ^{\alpha
(x)}dx-C(\sigma )\int_{\mathbb{R}^{N}}\mu w(x)\left\vert u\right\vert
^{\gamma (x)}dx \\
&\geq &2c_{1}\left\Vert u\right\Vert ^{\alpha ^{-}}-\left\vert u\right\vert
_{L_{a}^{\delta (\cdot )}}^{\delta (\xi _{1})}-\sigma \left\vert
u\right\vert _{L_{w}^{\alpha (\cdot )}}^{\alpha (\xi _{2})}-C(\sigma
)\left\vert u\right\vert _{L_{w}^{\gamma (\cdot )}}^{\gamma (\xi _{3})}\text{
(where }\xi _{1},\xi _{2},\xi _{3}\in \mathbb{R}^{N}\text{)} \\
&\geq &2c_{1}\left\Vert u\right\Vert ^{\alpha ^{-}}-\left\vert u\right\vert
_{L_{a}^{\delta (\cdot )}}^{\delta ^{+}}-\sigma \left\vert u\right\vert
_{L_{w}^{\alpha (\cdot )}}^{\alpha ^{+}}-C(\sigma )\left\vert u\right\vert
_{L_{w}^{\gamma (\cdot )}}^{\gamma ^{+}}-C_{1} \\
&\geq &2c_{1}\left\Vert u\right\Vert ^{\alpha ^{-}}-\beta _{k}^{\delta
^{+}}\left\Vert u\right\Vert ^{\delta ^{+}}-\sigma \beta _{k}^{\alpha
^{+}}\left\Vert u\right\Vert ^{\alpha ^{+}}-C(\sigma )\beta _{k}^{\gamma
^{+}}\left\Vert u\right\Vert ^{\gamma ^{+}}-C_{1} \\
&\geq &c_{1}\left\Vert u\right\Vert ^{\alpha ^{-}}-C_{2}\beta _{k}^{\delta
^{+}}\left\Vert u\right\Vert ^{\gamma ^{+}}-C_{3} \\
&=&c_{1}\left(2C_{2}\frac{1}{c_{1}}\beta _{k}^{\delta ^{+}}\right)^{\alpha ^{-}/(\alpha
^{-}-\gamma ^{+})}-C_{2}\beta _{k}^{\delta ^{+}}(2C_{2}\frac{1}{c_{1}}\beta
_{k}^{\delta ^{+}})^{\gamma ^{+}/(\alpha ^{-}-\gamma ^{+})}-C_{3} \\
&=&\frac{c_{1}}{2}\left(2C_{2}\frac{1}{c_{1}}\beta _{k}^{\gamma ^{+}}\right)^{\alpha
^{-}/(\alpha ^{-}-\gamma ^{+})}-C_{3}\rightarrow +\infty\ \text{ (as }%
k\rightarrow \infty \text{),}
\end{eqnarray*}%
because $1<\delta ^{+}<\alpha ^{-}<\gamma ^{+}$ and $\beta _{k}\rightarrow
0^{+}$ as $k\rightarrow \infty $. Therefore, $b_{k}\rightarrow +\infty $ as $%
k\rightarrow \infty $.

Now, we show that $(A_{2})$ holds. By ($\mathcal{H}_{f}^{2}$) and (\ref{w22}%
), we deduce that
\begin{equation*}
\Phi (tu)\rightarrow -\infty \text{ as }t\rightarrow +\infty ,\quad \forall
h\in V_{k}^{-}\text{ with}\parallel u\parallel =1,
\end{equation*}%
which implies that $(A_{2})$ holds.

We conclude that the proof of Theorem 1.2 is complete. $\square $

\subsection{Proof of Theorem 1.3}

\textbf{Proof of Theorem 1.3 } We first prove that $\Phi $ is coercive
on $X$.

Note that $\mu =0$. Since $1<\delta ^{+}<\alpha ^{-}$, then%
\begin{equation*}
\Phi (u)\geq \frac{1}{q^{+}}\left\Vert u\right\Vert ^{\alpha
^{-}}-c\left\Vert u\right\Vert ^{\delta ^{+}}\rightarrow +\infty \text{, as }%
\left\Vert u\right\Vert \rightarrow \infty .
\end{equation*}

By Lemma 4.4, $\Phi $ is weakly lower semi-continuous. Then $\Phi $  attains
its minimum on $X$, which provides a solution of problem $(\mathcal{E})$.

Since $\Phi $ is coercive, then $\Phi $ satisfies (PS) condition on $X$.
From assumption ($\mathcal{A}_{1}$)-(i), $\Phi $ is an even functional.
Denote by $\gamma (A)$ the genus of $A$ (see [\cite{2}], p. 215). Set%
\begin{eqnarray*}
\Re  &=&\{A\subset X\backslash \{0\}\mid A\text{ is compact and }A=-A\}, \\
\Re _{k} &=&\{A\subset \Re \mid \gamma (A)\geq k\}, \\
c_{k} &=&\underset{A\in \Re _{k}}{\inf }\underset{u\in A}{\sup }\Phi
(u),\ k=1,2,...
\end{eqnarray*}%
We have
\begin{equation*}
-\infty <c_{1}\leq c_{2}\leq \cdots \leq c_{k}\leq c_{k+1}\leq \cdots .
\end{equation*}

We prove in what follows that $c_{k}<0$ for every $k$.

 For fixed $k$, we can choose a $k$-dimensional linear subspace $E_{k}$
of $X$ such that $E_{k}\subset C_{0}^{\infty }(B_{R})$. As the norms on $%
E_{k}$ are equivalent, for any given $\delta _{0}>0$, there
exists $\rho _{k}\in (0,1)$ such that $u\in E_{k}$ with $\left\Vert
u\right\Vert \leq $ $\rho _{k}$ implies $\left\vert u\right\vert _{L^{\infty
}}\leq \delta _{0}$. Set%
\begin{equation*}
S_{\rho _{k}}^{(k)}=\{u\in E_{k}\mid \left\Vert u\right\Vert =\rho _{k}\}%
\text{.}
\end{equation*}

From the compactness of $S_{\rho _{k}}^{(k)}$, there exists $\theta
_{k}>0$ such that
\begin{equation*}
\int_{%
\mathbb{R}
^{N}}F(x,u)dx=\int_{%
\mathbb{R}
^{N}}\frac{\lambda a(x)\left\vert u\right\vert ^{\delta (x)}}{\delta (x)}%
dx\geq \theta _{k}\text{, }\quad\forall u\in S_{\rho _{k}}^{(k)}\text{.}
\end{equation*}

For $u\in S_{\rho _{k}}^{(k)}$ and $t\in (0,1)$, we have
\begin{eqnarray*}
\Phi \left( tu\right) &=&\int_{%
\mathbb{R}
^{N}}\mathscr{A}(x,\nabla tu)dx+\int_{%
\mathbb{R}
^{N}}\frac{1}{\alpha (x)}\left\vert tu\right\vert ^{\alpha (x)}dx-\int_{%
\mathbb{R}
^{N}}\frac{\lambda a(x)\left\vert tu\right\vert ^{\delta (x)}}{\delta (x)}dx
\\
&\leq &C_{1}\frac{t^{\alpha ^{-}}}{\alpha ^{-}}\rho _{k}^{\alpha
^{-}}-t^{\delta ^{+}}\theta _{k}.
\end{eqnarray*}

As $1<\delta ^{+}<\alpha ^{-}$, we can find $t_{k}\in (0,1)$ and $%
\varepsilon _{k}>0$ such that%
\begin{equation*}
\Phi \left( t_{k}u\right) \leq -\varepsilon _{k}<0,\text{ }\quad\forall u\in
S_{\rho _{k}}^{(k)},
\end{equation*}%
that is%
\begin{equation*}
\Phi \left( u\right) \leq -\varepsilon _{k}<0,\text{ }\quad\forall u\in
S_{t_{k}\rho _{k}}^{(k)}.
\end{equation*}

Obviously, $\gamma (S_{t_{k}\rho _{k}}^{(k)})=k$, so $c_{k}\leq -\varepsilon
_{k}<0.$

By the genus theory (see [\cite{2}], p. 219, Theorem 3.3), each $c_{k}$ is
a critical value of $\Phi $, hence there is a sequence of solutions \{$\pm
u_{k}\mid k=1,2,\cdots $\} such that $\Phi (\pm u_{k})<0$.

It only remains to prove that $c_{k}\rightarrow 0^{-}$ as $k\rightarrow \infty $.

Since $\Phi $ is coercive, there exists a constant $R>0$ such that $\Phi
(u)>0$ when $\left\Vert u\right\Vert \geq R$. Taking arbitrarily $A\in \Re
_{k}$, we have $\gamma (A)\geq k$. Let $Y_{k}$ and $Z_{k}$ be the subspaces of $%
X$ as mentioned in (\ref{8.1}). According to the properties of genus we know
that $A\cap Z_{k}$ $\neq \varnothing $.

Let $\beta _{k}=\sup \left\{ \left\vert \Phi _{f}(u)\right\vert \right\vert\
u\in Z_{k},\,\left\Vert u\right\Vert \leq R\}$. By Lemma 5.1, we have $\beta
_{k}\rightarrow 0$ as $k\rightarrow \infty $. For all $u\in Z_{k}$ with $%
\left\Vert u\right\Vert \leq R,$ we have
\begin{equation*}
\Phi (u)=\Phi _{\mathscr{A}}(u)+\Phi _{\alpha }(u)-\Phi _{f}(u)\geq -\Phi
_{f}(u)\geq -\beta _{k}\text{.}
\end{equation*}

Hence $\underset{u\in A}{\sup }\Phi (u)\geq -\beta _{k}$, and then $%
c_{k}\geq -\beta _{k}$. We conclude that $c_{k}\rightarrow 0^{-}$ as $%
k\rightarrow \infty $. $\square $

\section{Perspectives and open problems}
We now address to the readers several comments, perspectives, and open problems.

\smallskip
(i) Hypothesis ($\mathcal{A}_{1}$) (iv) establishes that problem $(\mathcal{E})$ is described in a {\it subcritical} setting. To the best of our knowledge, there is no result in the literature corresponding to the following {\it almost critical} framework described in what follows. Assume that condition
$q(\cdot )\ll \min \{N,p^{\ast }(\cdot )\}$ in ($\mathcal{A}_{1}$) (iv) is replaced with the following hypothesis: there exists a finite set $A\subset {\mathbb R}^N$ such that $q(a)=\min \{N,p^{\ast }(a )\}$ for all $a\in A$ and
$q(x )< \min \{N,p^{\ast }(x )\}$ for all $x\in{\mathbb R}^N\setminus A$.

\smallskip
{\it Open problem.} Study if
 Theorems 1.1--1.3 established in this paper still remain true in the above almost critical abstract setting.

\smallskip
(ii)  Another very interesting research direction is to extend the approach developed in this paper to the case of {\it double phase} problems studied by Mingione {\it et al.} [\cite{mingi1, mingi0, mingi2}]. This corresponds to the following non-homogeneous potential
$$
\mathscr{A}(x,\xi )=\frac{a(x)}{p(x)}\,|\xi|^{p(x)}+\frac{b(x)}{q(x)}\,|\xi|^{q(x)},$$
where the coefficients $a(x)$ and $b(x)$ are non-negative and at least one is strictly positive for all $x\in{\mathbb R}^N$. At this stage, we do not know any multiplicity results for double phase problems of this type.

We also refer to the pioneering papers by Marcellini [\cite{marce1, marce2}] on $(p,q)$-growth conditions, which involve integral functionals of the type
$$W^{1,1}\ni u\mapsto\int_\Omega f(x,\nabla u)dx,$$
where $\Omega\subseteq{\mathbb R}^N$ is an open set.
The integrand $f:\Omega\times{\mathbb R}^N\rightarrow {\mathbb R}$ satisfied unbalanced polynomial growth conditions of the type
$$|\xi|^p \lesssim f(x,\xi)\lesssim |\xi|^q+1\quad\mbox{with $1<p<q$,}$$
for every $x\in\Omega$ and $\xi\in{\mathbb R}^N$.

\smallskip
(iii) The differential operator $\mathscr{A}(x,\xi )$ considered in problem $(\mathcal{E})$ falls in the realm of those related to the so-called Musielak-Orlicz spaces (see [\cite{19a, orli}]), more in general, of the operators having non-standard growth conditions (which are widely considered in the calculus of variations). These function spaces are Orlicz spaces whose defining Young function exhibits an additional dependence on the $x$ variable. Indeed, classical Orlicz spaces $L^\Phi$ are defined requiring that a member function $f$ satisfies
$$\int_\Omega \Phi(|f|)dx<\infty,$$
where $\Phi(t)$ is a Young function (convex, non-decreasing, $\Phi(0)=0$). In the new case of Musielak-Orlicz spaces, the above condition becomes
$$\int_\Omega \Phi(x,|f|)dx<\infty.$$

The problems considered in this paper are indeed driven by the function
\begin{equation}\label{dphase}\Phi(x,|\xi|):=\left\{
\begin{array}{lll}
& |\xi|^{p(x)}&\quad\mbox{if}\ |\xi|\leq 1\\
& |\xi|^{q(x)}&\quad\mbox{if}\ |\xi|\geq 1.
\end{array}\right.
\end{equation}
When $p(x)=q(x)$ we find the so-called variable exponent spaces, which are defined by
$$\Phi(x,|\xi|):=|\xi|^{p(x)}.$$

We conclude these comments by saying that the present paper is concerned with a double phase variant of the operators stemming from the energy generated by the function defined in (\ref{dphase}).

\smallskip
(iv) An interesting double phase type operator  considered in the papers of Baroni, Colombo and Mingione [\cite{mingi1,mingi0,mingi2}], addresses functionals of the type
\begin{equation}\label{eenerg}w\mapsto \int_\Omega (|\nabla w|^p+a(x)|\nabla w|^q)dx,\end{equation}
where $a(x)\geq 0$. The meaning of this functional is also to give a sharper version of the following energy
$$w\mapsto \int_\Omega |\nabla w|^{p(x)}dx,$$
thereby describing sharper phase transitions. Composite materials with locally different hardening exponents $p$ and $q$ can be described using the energy defined in (\ref{eenerg}). Problems of this type are also motivated by applications to elasticity, homogenization, modelling of strongly anisotropic materials, Lavrentiev phenomenon, etc.

Accordingly, a new double phase model can be given by
$$\Phi_d(x,|\xi|):=\left\{
\begin{array}{lll}
& |\xi|^{p}+a(x)|\xi|^{q}&\quad\mbox{if}\ |\xi|\leq 1\\
& |\xi|^{p_1}+a(x)|\xi|^{q_1}&\quad\mbox{if}\ |\xi|\geq 1,
\end{array}\right.
$$
with $a(x)\geq 0$.

\medskip
{\bf Acknowledgments.} This research was supported by the Slovenian Research Agency grants
P1-0292, J1-8131, N1-0064, N1-0083, and N1-0114.

\end{document}